\input amstex.tex
\input amsppt.sty
\magnification=\magstep 1
\pageheight{46pc}
\pagewidth{30pc}
\topmatter
\title  Maximum principles for a class of nonlinear 
second order elliptic differential equations \endtitle
\rightheadtext { Maximum principles and overdetermined problems }
\author G. Porru, A. Tewodros and S. Vernier-Piro \endauthor
\leftheadtext {G. Porru, A. Tewodros and S. Vernier-Piro}
\address Giovanni Porru and Stella Vernier-Piro: Dipartimento di Matematica, 
Via Ospedale 72, 09124, Cagliari, Italy. \endaddress
\address Amdeberhan Tewodros: Mathematics Department, Temple University, 
Philadelphia, PA 19122, USA. \endaddress
\abstract In this paper we investigate maximum principles
for functionals defined on solutions to special partial differential
equations of elliptic type, extending results by Payne and Philippin.
We apply such maximum principles to investigate one overdetermined problem.
\endabstract
\endtopmatter
\document
\subhead 1. Introduction \endsubhead
\medskip
We consider classical  solutions  $u=u(x)$  of
the quasilinear second order equation
$$
\bigl(g(q^2)u_i\bigr)_i=h(q^2)                               \leqno(1.1)$$
in domains  $\Omega\subset R^N$.  Here and in the sequel the
subindex $i\ (i=1,...,N)$ denotes partial differentiation with  respect  to $x^i$,
the summation convention (from $1$ to $N$) over repeated indices is in effect,
$q^2=u_iu_i$, $g$ and $h$ are  two  smooth  functions. In order
for equation (1.1) to be elliptic we suppose $g>0$ and $G>0,$ where
$$
G(\xi)=g(\xi)+2\xi g'(\xi).                                 \leqno(1.2)$$
Following Payne and Philippin [5,6,7]  we  derive some maximum  principles for
functionals $\Phi(u,q)$ defined  on  solutions  $u$  of  equation  (1.1). Of
course, there are infinitely many choices for such functionals.  In  order to
exploit the corresponding maximum principle for getting more  information on
$u$, not only the functional must  satisfy  a  maximum  principle,  but,  in
addition, there must exist some domain $\Omega$ and some solution $u$ of (1.1)
for which $\Phi (u,q)$ is a constant throughout $\Omega$.
Such maximum principles are named "best possible" maximum principles ([4]).
For applications of such "best possible" inequalities in fluid
mechanics, geometry and in other areas we refer to [4].
 
In [6] Payne and Philippin consider the functional
$$
\Phi(u,q)={1\over 2}\int_0^{q^2}{G(\xi)\over h(\xi)}d\xi-u       \leqno(1.3)$$
and prove that, if $u(x)$ satisfies the equation (1.1) then $\Phi(u,q)$
assumes its maximum value either on the boundary of $\Omega$ or when $q=0$.
If $u=u(x^1)$ is a function depending on one variable only and if it is a
solution of (1.1) then the corresponding $\Phi(u,q)$ is a constant.
In the same paper [6], Payne and Philippin define
$$
\Psi(u,q)={N\over 2}\int_0^{q^2}{G(\xi)\over h(\xi)}d\xi-u       \leqno(1.4)$$
and prove that, if $\Psi(u,q)$ is computed on any solution of equation
(1.1) then it assumes its maximum value on the boundary of $\Omega$.
In case $h$ is a nonvanishing constant and $u=u(r)$ is the radial solution
of equation (1.1) satisfying $u'(0)=0$,
then $\Psi(u,q)$ is a constant.
These results have been extended to more general equations in [7].
 
In Section 2 of this paper we exhibit a new class of functionals
which satisfy "best possible" maximum principles.
These functionals are expressed in
terms of solutions to an ordinary differential equation related to (1.1).
 
In Section 3 we consider the equation:
$$
\bigl(g(q^2)u_i\bigr)_i=N,\leqno(1.5)$$
where $g$ satisfies suitable hypotheses. Assume equation (1.5) has a smooth solution
$u(x)$ in a convex ringshaped domain $\Omega\subset R^N$ bounded externally
by a (hyper) surface $\Gamma_0$ and internally by a (hyper) surface $\Gamma_1$.
We show that, if such a solution satisfies the following (overdetermined)
boundary conditions
$$
u|_{\Gamma_0}=0, \; \; \; \; u|_{\Gamma_1}=-c_1,$$
$$
u_n|_{\Gamma_0}=q_0, \; \; \; \; u_n|_{\Gamma_1}=0,$$
where $c_1$ and $q_0$ are positive free constants, then ${\Gamma_0}$ and
${\Gamma_1}$ must be two concentric N-spheres. Similar problems have been
investigated by several authors. In [9] Philippin and Payne discussed the
equation:
$$
\bigl(q^{N-2}u_i\bigr)_i=0$$
under the boundary conditions
$$
u|_{\Gamma_0}=0, \; \; \; \; u|_{\Gamma_1}=-1,                  \leqno(1.6)$$
$$
u_n|_{\Gamma_0}=q_0, \; \; \; \; u_n|_{\Gamma_1}=-q_1,           \leqno(1.7)$$
where $q_0$ and $q_1$ are free constants. They proved
that if this problem is solvable then
$\Omega$ must be radially symmetric. In [8] Philippin solved
tha same problem in case the equation is $\Delta u=0$ and
the boundary conditions are (1.6),
(1.7). In [10] Porru and Ragnedda investigated the above problem when the
equation is
$$
\bigl(q^{p-2}u_i\bigr)_i=0,$$
$p>1$, again under conditions (1.6), (1.7).
The case when $\Omega$ is a bounded simple
connected domain has been studied by Serrin in [12]. By using the moving
plane method he has found that if
$u(x)$ is a smooth solution of equation (1.1) and satisfies
$$
u|_{\partial\Omega}=0, \; \; \; \; u_n|_{\partial\Omega}=q_0,     \leqno(1.8)$$
($q_0\not=0$) then $\Omega$ must be a sphere. The same result has been found by
Weinberger [13] for the special case $\Delta u=-1$ by using a different
method. Extending Weinberger's method, Garofalo and Lewis [1] have
solved the overdetermined problem (1.5), (1.8) allowing $u(x)$ to be a
generalized solution.
\medskip
\subhead 2. Maximum principles \endsubhead          

Let us prove first some preliminary lemmas.

\proclaim {Lemma 2.1} Let $z(t)$ be either a strictly convex or a strictly
concave $C^2$ function in $(t_0,t_1)$ and let $\psi (s)$ be the inverse
of $z'(t)$. Then the function of $t$
$$
\int_{z'(t_0)}^{z'(t)}s\psi'(s)\, ds-z(t)$$
is a constant on $(t_0,t_1)$. \endproclaim
 
\demo {Proof} The proof is trivial.
If we replace $s$ by $z'(\tau)$ in the above integral we obtain
$$
\int_{z'(t_0)}^{z'(t)}s\psi'(s)\, ds-z(t)=\int_{t_0}^{t}z'(\tau)\,
d\tau-z(t)=-z(t_0).$$
The lemma is proved.  \enddemo
 
\proclaim {Lemma 2.2} Let $u=u(x)$ be a smooth function 
satisfying $\nabla u\not=0$ in
$\Omega\subset R^N$, and let $\psi =\psi (s)$ be a smooth function in
$(0,\infty)$. Then we have in $\Omega$
$$
\psi^2u_{ih}u_{ih}                                           \leqno(2.1)$$
$$
\geq\bigl(\psi^2-(\psi')^2q^2\bigr)q_iq_i+
\bigl(2\psi'q-2\psi\bigr)q_iu_i-Nq^2+2\psi q\Delta u,$$
where $q=|\nabla u|$ and $\psi=\psi(q)$. Furthermore, equality holds in $(2.1)$
throughout $\Omega$ if and only if $\psi^2(q)=|x-x_0|^2$.
\endproclaim
 
\demo {Proof}
Let $\delta^{ih}$ be the Kronecker delta. By
$$
\sum_{i,h}^{1,N}\Bigl(\bigl({\psi(q)\over q}u_i\bigr)_h-\delta^{ih}\Bigr)^2
\geq 0                                                          \leqno(2.2)$$
it follows
$$
\sum_{i,h}^{1,N}\Bigl({{\psi}\over q}u_{ih}+{\psi'\over q}q_hu_i-{{\psi}
\over{q^2}}q_hu_i-\delta^{ih}\Bigr)^2\geq 0.$$
Easy computations give
$$
{{\psi^2}\over{q^2}}u_{ih}u_{ih}-\bigl({{\psi^2}\over{q^2}}
-(\psi')^2\bigr)q_iq_i-\bigl(2{\psi'\over q}-2{\psi
\over{q^2}}\bigr)q_iu_i+N-{{2\psi}\over q}\Delta u\geq0,$$
where the  identities  $u_iu_i=q^2$,  $u_{ih}u_h=qq_i$
have been used. Inequality (2.1) follows. Of course, we have equality
in (2.1) if and only if equality holds in (2.2), that is, if and only if 
$$
\bigl({\psi(q)\over q}u_i\bigr)_h=
\delta^{ih},\ \ i,h=1,\ldots,N.$$
By solving the last equations we obtain
$$
{\psi(q)\over q}u_i=x^i-x^i_0, \ \ i=1,\ldots,N.               \leqno(2.3)$$
Since $u_iu_i=q^2$, these equalities imply
$$
\psi^2(q)=|x-x_0|^2. $$
The lemma has been proved. \enddemo
 
\proclaim {Corollary} If we have equality in $(2.1)$ and 
if $\psi(s)$ is either
strictly increasing or strictly decreasing then $u(x)$ must be
a radial function. \endproclaim
 
\demo {Proof} Since $\psi^2(q)=|x-x_0|^2$, by (2.3) we find:
$$
\nabla u=H(|x-x_0|)(x-x_0),$$
where $H$ is an appropriate function of one variable only. The result follows.
\enddemo

Let us come to equation (1.1). The function $g$ is assumed to be smooth on
$(0,\infty)$ and to satisfy
$$
g(\xi)>0,\ \ \lim_{\xi\to 0}g(\xi^2)\xi=0,\ \ \lim_{\xi\to \infty}g(\xi^2)\xi=
\infty,\ \ G(\xi)>0,                                           \leqno(2.4)$$
where $G(\xi)$ is defined as in (1.2). The function $h$ is supposed to be
smooth in $[0,\infty)$ and to satisfy
$$
h(\xi)>0,\ \ h'(\xi)\geq 0\ \ \forall\xi\geq 0.                \leqno(2.5)$$
The case $h(\xi)<0,\  h'(\xi)\leq 0$
can be reduced to the case in above by changing $u$ with $-u$ in (1.1).
 
Define the ordinary differential equation
$$
\Bigl(t^{N-1}g(\phi^2)\phi\Bigr)'=t^{N-1}h(\phi^2).             \leqno(2.6)$$
Observe that, if $u(r)$, $r=|x|$, is a radial solution of equation (1.1)
for $r_1<r<r_2$ then $\phi(t)=u'(t)$ is a solution of (2.6) for $r_1<t<r_2$.
 
\proclaim {Lemma 2.3} Assume conditions $(2.4),\ (2.5)$. 
Given $t_0\geq 0$, let
$\phi(t)$ be the solution of equation $(2.6)$ satisfying $\phi(t_0)=0$. If
$(t_0,t_1)$ is the maximal interval of existence for $\phi(t)$ then
$\phi'(t)>0$ on $(t_0,t_1)$ and $\phi(t)\to \infty$ as $t\to t_1$.
Here $t_1$ may be finite or $\infty$. \endproclaim
 
\demo {Proof} This lemma is probably known, but we give a proof for completeness.
From the equation (2.6) and the condition $\phi(t_0)=0$ one finds
$\phi'(t)>0$ on the (maximal) interval $(t_0,a)$,
with $a\leq t_1$. We claim that $a=t_1$. By contradiction, let $a<t_1$, so
$\phi'(a)=0$. Since $h$ is nondecreasing, the function $h(\phi^2(t))$
is nondecreasing on $(t_0,a)$. Hence, integration of (2.6) on $(t_0,t)$,
$t\leq a$, yields
$$
t^{N-1}g(\phi^2)\phi \leq {t^N-t_0^N\over N}h(\phi^2)\leq{t^N\over
N}h(\phi^2).                                                        \leqno(2.7)$$
Insertion of (2.7) into (2.6) rewritten as $G(\phi^2)\phi'+{N-1\over t}
g(\phi^2)\phi=h(\phi^2)$ leads to
$$
G(\phi^2)\phi'\geq{1\over N}h(\phi^2).                                 \leqno(2.8)$$
At $t=a$, (2.8) implies $\phi'(a)>0$, which contradicts the assumption
$\phi'(a)=0$. Hence $\phi(t)$ is strictly increasing on $(t_0,t_1)$.
 
If $t_1$ is finite then
$\phi(t)\to \infty$ as $t\to t_1$ because of the maximality of the interval
$(t_0,t_1)$. Let $t_1=\infty.$  For $t\geq t_0$, (2.6) implies
$$
\Bigl(t^{N-1}g(\phi^2)\phi\Bigr)'\geq t^{N-1}h(0).$$
Integrating over $(t_0,t)$ we find
$$
g(\phi^2)\phi\geq t\Bigl(1-\Bigl({t_0\over t}\Bigr)^N\Bigr){h(0)\over N}.$$
Taking into account conditions (2.4), the above inequality implies that
$\phi(t)\to \infty$ as $t\to \infty$.
The lemma is proved.
\enddemo

Observe that equation (1.1) may be rewritten as
$$
\Delta u+f(q)u_iu_ju_{ij}=k(q),                                \leqno (2.9)$$
where $f(q)=2g'(q^2)/g(q^2)$ and $k(q)=h(q^2)/g(q^2)$.
The ordinary differential equation (2.6) in terms of $f$ and $k$ reads as
$$
\phi'\bigl(1+\phi^2f(\phi)\bigr)+{N-1\over t}\phi=k(\phi).    \leqno(2.10)$$
 
\proclaim {Theorem 2.1} Assume conditions $(2.4),\ (2.5)$.
Given $t_0\geq 0$, let $\phi(t)$
be the solution of equation $(2.10)$ satisfying $\phi(t_0)=0$, and let $\psi(s)$
be the inverse function of $\phi(t)$. If $u(x)$ is a solution of equation
$(2.9)$ such that $\nabla u \not=0$ in $\Omega$
then the function
$$
\Phi (u,q)=\int_0^{q(x)}s\psi'(s)\, ds-u(x)                      \leqno(2.11)$$
assumes its maximum value on the boundary of $\Omega$.
Moreover, $\Phi (u,q)$ is a constant if $u(x)$ is the
(radially symmetric) solution $u(x)=F(|x|)$, $F'(t)=\phi(t)$.
\endproclaim
 
\demo {Proof} By Lemma 2.3, $\phi(t)$ is strictly increasing, hence
the second part of the theorem follows by Lemma 2.1 when $z'(t)=\phi(t)$.
For proving the first part we put
$$
v(x)=\int_0^{q(x)}s\psi'(s)\, ds-u(x),$$
where $u(x)$ is a solution of equation (2.9). We have
$$
v_i=\psi'(q)qq_i-u_i,\ \ i=1,\ldots,N.                           \leqno(2.12)$$
By (2.12) we obtain
$$
v_{ij}=\psi'(q_iq_j+qq_{ij})+\psi''qq_iq_j-u_{ij},\ \ i,j=1,\ldots,N.$$
From the identities $qq_i=u_{ih}u_h$ we get
$$
q_iq_j+qq_{ij}=u_{ih}u_{jh}+u_{ijh}u_h.$$
Consequently, we find
$$
v_{ij}=\psi'(u_{ih}u_{jh}+u_{ijh}u_h)+\psi''qq_iq_j-u_{ij}.      \leqno(2.13)$$
Let us define
$$
a^{ij}=\delta^{ij}+f(q)u_iu_j,\ \ i,j=1,\ldots,N,                \leqno(2.14)$$
where $\delta^{ij}$ is the Kronecker delta. In virtue of conditions (2.4)
the matrix $[a^{ij}]$ is positive definite.
By using (2.13) and (2.14) we find
$$
{1\over{\psi'}}a^{ij}v_{ij}=u_{ih}u_{ih}+fq^2q_iq_i$$
$$
+a^{ij}u_{ijh}u_h+{{\psi''}\over
{\psi'}}q(q_iq_i+fq_iu_iq_ju_j)-{k\over{\psi'}},$$
where the equation (2.9) rewritten as $a^{ij}u_{ij}=k$ has been used.
Easy computations yield
$$
a^{ij}u_{ijh}=(a^{ij}u_{ij})_h-(a^{ij})_hu_{ij}=k'q_h-
f'q_hqq_iu_i-2fu_{ih}qq_i,$$
where $f'=f'(q)$ and $k'=k'(q)$. Hence
$$
{1\over{\psi'}}a^{ij}v_{ij}                           \leqno(2.15)$$
$$
=u_{ih}u_{ih}-fq^2q_iq_i+k'q_iu_i-
f'qq_iu_iq_ju_j+{{\psi''}\over{\psi'}}q(q_iq_i+fq_iu_iq_ju_j)-{k\over{\psi'}}.$$
Equality (2.15) and inequality (2.1) give
$$
{1\over{\psi'}}a^{ij}v_{ij}                                   \leqno(2.16)$$
$$
\geq\Bigl(1-{(\psi')^2\over{\psi^2}}q^2\Bigr)q_iq_i+\Bigl(2{\psi'\over{\psi^2}}q-
{2\over\psi}\Bigr)q_iu_i-{N\over{\psi^2}}q^2+{2\over\psi}kq-
{2\over\psi}fq^2q_iu_i$$
$$
-fq^2q_iq_i+k'q_iu_i-f'qq_iu_iq_ju_j+
{{\psi''}\over{\psi'}}q(q_iq_i+fq_iu_iq_ju_j)-{k\over{\psi'}},$$
where the equation $\Delta u=k-fqq_iu_i$ has been used. By (2.12) we obtain
$$
q_ju_j={v_ju_j\over{\psi'q}}+{q\over{\psi'}},\,\,\,\,\, q_iq_i=\Bigl(q_i+
{u_i\over{\psi'q}}\Bigr){v_i\over{\psi'q}}+{1\over{(\psi')^2}}.\leqno(2.17)$$
Insertion of (2.17) into (2.16) yields:
$$
{1\over{\psi'}}a^{ij}v_{ij}                       \leqno(2.18)$$
$$
\geq{1\over{(\psi')^2}}+{q^2\over{\psi^2}}(1-N)-{2q\over{\psi\psi'}}(1+fq^2)
+{2kq\over\psi}+{k'q\over{\psi'}}-$$
$$
-{f'q^3\over{(\psi')^2}}-{fq^2\over{(\psi')^2}}+{\psi''\over{(\psi')^3}}
q(1+fq^2)-{k\over{\psi'}}-{1\over{\psi'}}b^iv_i,$$
where $b^i$ (i=1,$\ldots$,N) is a regular vector field (recall that, by
assumption, $q(x)>0$).
From equation (2.10) with $t=\psi (q)$ (and, consequently, $\phi(t)=q$) we find
$$
k={1\over{\psi'}}(1+fq^2)+{{N-1}\over{\psi}}q.                   \leqno(2.19)$$
By using equation (2.19), inequality (2.18) becomes:
$$
{1\over{\psi'}}a^{ij}v_{ij}+{1\over{\psi'}}b^iv_i            \leqno(2.20)$$
$$
\geq{\psi''\over{(\psi')^3}}q(1+fq^2)-{2fq^2\over{(\psi')^2}}-
{f'q^3\over{(\psi')^2}}+{k'q\over{\psi'}}-{N-1\over{\psi\psi'}}q
+{N-1\over{\psi^2}}q^2.$$
Differentiation with respect to $t$ in equation (2.10) yields:
$$
\phi''\bigl(1+\phi^2f\bigr)+2(\phi')^2\phi f+(\phi')^2\phi^2f'
+{N-1\over t}\phi'-{N-1\over{t^2}}\phi=k'\phi'.                  \leqno(2.21)$$
Since $\psi(q)$ is the inverse of $\phi(t)$, we have
$$
\phi'(t)={1\over{\psi'(q)}},\,\,\,\phi''(t)=-{{\psi''(q)}\over{(\psi'(q))
^3}}.$$
Hence, the equation (2.21) may be rewritten as
$$
{k'\over{\psi'}}=
-{{\psi''}\over{(\psi')^3}}(1+fq^2)+{2fq\over{(\psi')^2}}+
{f'q^2\over{(\psi')^2}}+{{N-1}\over{\psi\psi'}}-{{N-1}\over{\psi^2}}q.
                                                                \leqno(2.22)$$
Insertion of (2.22) into (2.20) yields
$$
a^{ij}v_{ij}+b^iv_i\geq 0.                                   \leqno(2.23)$$
The theorem follows by (2.23) and the classical maximum principle
[11,2].
\enddemo
 
\it Remark 2.1 \rm If equality holds in (2.23) then equality
holds in (2.1) and, by the corollary to Lemma 2.2, $u(x)$
is a radial function.
\smallskip
\it Remark 2.2 \rm If $h$ is a (positive) constant and if $t_0=0$ then the
functional (2.11) is the same as that defined in (1.4).
\medskip
In case of dimension two, Theorem 2.1 can be improved. In fact, we have the
following
\proclaim {Theorem 2.2} Under the same assumptions as in Theorem 2.1, if N=2
the function $v(x)=\Phi (u,q)$ defined in $(2.11)$
assumes its maximum value and its minimum value on the boundary of $\Omega$.
\endproclaim
 
\demo {Proof} The proof of this theorem is the same as that of Theorem 2.1 up
to equation (2.15). At this point, instead of using inequality (2.1),
we make use of the following equality (true for $N=2$ only, ([7]) p. 43))
$$
u_{ih}u_{ih}=(\Delta u)^2+2q_iq_i-2\Delta uq_i{u_i\over q}.\leqno(2.24)$$
Since $\Delta u=k-fqq_iu_i$, equality (2.24) yields
$$
u_{ih}u_{ih}=\bigl(k-fqq_iu_i\bigr)^2+2q_iq_i-2kq_i{u_i\over q}+
2fq_iu_iq_ju_j.                                                \leqno(2.25)$$
Insertion of (2.25) into (2.15) and use of equations (2.17) lead to
$$
{1\over{\psi'}}a^{ij}v_{ij}                                  \leqno(2.26)$$
$$
=\bigl(k-{fq^2\over{\psi'}}\bigr)^2+{2\over{(\psi')^2}}-{3k\over{\psi'}}+
{fq^2\over{(\psi')^2}}$$
$$
+{k'q\over{\psi'}}-{f'q^3\over{(\psi')^2}}+{\psi''q\over{(\psi')^3}}(1+fq^2)
-{1\over{\psi'}}d^iv_i,$$
where
$d^i$ (i=1, 2) is a regular vector field. By equation (2.19) with $N=2$ we have
$$
k={1\over{\psi'}}+{1\over{\psi'}}fq^2+{1\over{\psi}}q.         \leqno(2.27)$$
By (2.22) we have
$$
{k'\over{\psi'}}=-{\psi''\over{(\psi')^3}}(1+fq^2)+{2fq\over{(\psi')^2}}+
{f'q^2\over{(\psi')^2}}+{1\over{\psi\psi'}}-{q\over{\psi^2}}.  \leqno(2.28)$$
Insertion of (2.27) and (2.28) into (2.26) leads to
$$
a^{ij}v_{ij}+d^iv_i=0.$$
The theorem follows by classical maximum principles [11,2].
\enddemo

\subhead 3. An overdetermined problem \endsubhead

Throughout this section we assume $\Omega\subset R^N$ to be a smooth
ringshaped domain bounded externally by a (hyper) surface
$\Gamma_0$ and internally by a (hyper) surface $\Gamma_1$. We also
suppose $\Gamma_0$ and $\Gamma_1$ enclose convex domains $\Omega_0$ and
$\Omega_1$, respectively. In $\Omega$ we investigate the following
(overdetermined) problem:
$$
\bigl(g(q^2)u_i\bigr)_i=N,                                        \leqno(3.1)$$
$$
u|_{\Gamma_0}=0, \; \; \; \; u|_{\Gamma_1}=-c_1,                  \leqno(3.2)$$
$$
u_n|_{\Gamma_0}=q_0, \; \; \; \; u_n|_{\Gamma_1}=0,               \leqno(3.3)$$
where $c_1$ and $q_0$ are two positive free constants and the subindex $n$
denotes normal external differentiation. The function $g(\xi)$ is assumed to
satisfy conditions (2.4). According to (1.2) we have
$$
G(s^2)={d\over {ds}}\bigl(g(s^2)s\bigr).                         \leqno(3.4)$$
We are interested only in
smooth solutions of equation (3.1) whose gradient is nonvanishing in $\Omega$.

\proclaim {Theorem 3.1} If $(2.4)$ holds then problem $(3.1)$, $(3.2)$, 
$(3.3)$ is solvable if and only  if $q_0$, $c_1$ satisfy
$$
{1\over N}\int_0^{q_0}sG(s^2)ds < c_1 < \int_0^{q_0}sG(s^2)ds     \leqno(3.5)$$
and $\Gamma_0$, $\Gamma_1$ are two suitable concentric N-spheres.
\endproclaim
 
\demo {Proof} The ordinary differential equation (2.6) corresponding to
our partial differential equation (3.1) is
$$
\Bigl(t^{N-1}g(\phi^2)\phi\Bigr)'=\bigl(t^N\bigr)'.              \leqno (3.6)$$
Take $\alpha \geq 0$ and assume $\phi(\alpha)=0$. Integrating (3.6) over
$(\alpha,t)$ we obtain
$$
g(\phi^2)\phi=t(1-\alpha^Nt^{-N}).$$
Taking into account conditions (2.4) one concludes that $\phi(t)$ is
defined on $(\alpha,\infty)$.
If $\psi(s)$ denotes the inverse function of $\phi(t)$ then
by the last equation we have
$$
g(s^2)s=(1-\alpha^N\psi^{-N})\psi.                                 \leqno(3.7)$$
From (3.7) we get
$$
\psi'={G(s^2)\psi^N\over{\psi^N+(N-1)\alpha^N}}, \; \; \; \psi(0)=\alpha.
                                                                 \leqno(3.8)$$
 
Suppose problem (3.1), (3.2), (3.3) has a regular solution $u=u(x)$. Let us
consider the function $\psi(s)$ defined by (3.7) when $\alpha=0$. We find
$\psi(s)=g(s^2)s$ and $\psi'(s)=G(s^2)$. By Theorem 2.1, the function
$$
v(x)=\int_0^{q(x)}sG(s^2)ds-u(x)                                 \leqno(3.9)$$
attains its maximum value on $\Gamma_1\cup\Gamma_0$. On $\Gamma_1$
and on $\Gamma_0$ we have
$$
v_n=qG(q^2)q_n-u_n,                                            \leqno(3.10)$$
where the subindex $n$, as before, denotes normal external differentiation.
Equation (3.1) rewritten in normal coordinates reads as
$$
G(q^2)q_n+g(q^2)(N-1)Ku_n=N,$$
where $K$ is the mean curvature of the corresponding level surface. Since
$\Gamma_1$ is smooth and since $u_n=-q=0$ on $\Gamma_1$, we find
$$
G(q^2)q_n=N\ \ \ \hbox{\rm on}\  \Gamma_1.                      \leqno(3.11)$$
It follows that $v_n$ vanishes at each point in $\Gamma_1$. Hence, by Hopf's
second principle, $v(x)$ cannot take its maximum value on $\Gamma_1$
unless it is a constant in $\Omega$.
Consequently, such a maximum value is attained in $\Gamma_0$.
Using conditions (3.2), (3.3) and comparing the values of $v$ on $\Gamma_1$
and $\Gamma_0$ we find
$$
c_1\leq \int_0^{q_0}sG(s^2)ds.                               \leqno(3.12)$$
Using again conditions (3.2), (3.3) we obtain
$$
v(x)\leq\int_0^{q_0}sG(s^2)ds,$$
from which it follows
$$
\int_\Omega v(x)dx\leq(|\Omega_0|-|\Omega_1|)\int_0^{q_0}sG(s^2)ds.\leqno(3.13)$$
On the other hand, from equation (3.1) we find
$$
{1\over N}\Bigl(x^ju_jgu_i-gq^2x^i+(N-1)ugu_i\Bigr)_i={x^i\over N}\Bigl(
gqq_i-\bigl(gq^2\bigr)_i\Bigr)+x^iu_i+(N-1)u.                    \leqno(3.14)$$
By using Green's formula as well as conditions (3.2), (3.3) we obtain
$$
\int_\Omega\Bigl(x^ju_jgu_i-gq^2x^i+(N-1)ugu_i\Bigr)_idx=0.     \leqno(3.15)$$
Since
$$
gqq_i-\bigl(gq^2\bigr)_i=-(gq)_iq=-\Bigl(\int_0^qsG(s^2)ds\Bigr)_i$$
we find
$$
\int_\Omega{x^i\over N}\Bigl(gqq_i-\bigl(gq^2\bigr)_i\Bigr)dx=
-\int_\Omega{x^i\over N}\Bigl(\int_0^qsG(s^2)ds\Bigr)_idx.$$
By using again Green's formula, the boundary conditions (3.3)
and the well known equation
$$
\int_{\Gamma_0}x^in^ids=N|\Omega_0| $$
we obtain
$$
\int_\Omega{x^i\over N}\Bigl(gqq_i-\bigl(gq^2\bigr)_i\Bigr)dx=
-\int_0^{q_0}sG(s^2)ds|\Omega_0|+\int_\Omega\Bigl(\int_0^qsG(s^2)ds\Bigr)dx.
                                                                 \leqno(3.16)$$
By using once more Green's formula and the boundary conditions (3.2) we find
$$
\int_\Omega x^iu_idx=-c_1\int_{\Gamma_1}x^in^ids-N\int_\Omega u\,dx. $$
Since
$$
\int_{\Gamma_1}x^in^ids=-N|\Omega_1|, $$
the previous equation gives
$$
\int_\Omega x^iu_idx=c_1N|\Omega_1|-N\int_\Omega u\,dx.           \leqno(3.17)$$
Integration in (3.14) and use of (3.15), (3.16), (3.17) and (3.9) lead to
$$
\int_{\Omega}v(x)dx=\int_0^{q_0}sG(s^2)ds|\Omega_0|-c_1N|\Omega_1|.\leqno(3.18)$$
By (3.13) and (3.18) it follows
$$
{1\over N}\int_0^{q_0}sG(s^2)ds\leq c_1.                          \leqno(3.19)$$
If equality were to occur in (3.12) then equality would occur in
(3.13) and in (3.19), a contradiction. Hence we must have strict inequality in
(3.12) and in (3.19). Inequalities (3.5) have been proved.
 
Now let us consider the
function $\psi(s)$ defined by (3.7) for $\alpha >0$, and discuss the
following equation
$$
\int_0^{q_0}s\psi'(s)\,ds=c_1.                                  \leqno(3.20)$$
By (3.7) with $\alpha =0$ we find $s\psi'(s)=sG(s^2).$ Since $\psi(0)=\alpha$
and since $\psi(s)$ is increasing, for $s$ fixed, $\psi\to\infty$ as
$\alpha \to \infty$. As a consequence,
by (3.7) we infer that when $\alpha\rightarrow\infty$ then
${\alpha\over{\psi(s)}}\rightarrow 1$. Hence, by (3.8) it follows
that
$$
s\psi'(s)\rightarrow {1\over N}sG(s^2)\ \ \hbox{\rm as}\ \ \alpha\to \infty.$$
Therefore, since $s\psi'(s)$ decreases as $\alpha$ increases, and since
$q_0$ and $c_1$ satisfy inequalities (3.5), there is a unique positive
$\alpha$ which solves equation (3.20). By using this value of
$\alpha$ let us define
$$
v(x)=\int_0^{q(x)}s\psi'(s)\,ds-u(x).                           \leqno(3.21)$$
In virtue of conditions (3.2), (3.3) and (3.20),
the function $v(x)$ assumes
the same value on $\Gamma_0$ and on $\Gamma_1$. If $N=2$ then Theorem 2.2
implies that $v(x)$ is a constant in $\Omega$. For general $N$, let us
compute the normal derivative of $v$ on $\Gamma_1$. We find
$$
v_n=q\psi'(q)q_n-u_n.                                            \leqno(3.22)$$
By using (3.8) with $s=q$ we have
$$
\psi'q_n={\psi^N\over{\psi^N+(N-1)\alpha^N}}G(q^2)q_n.          \leqno(3.23)$$
By (3.23) and (3.11), recalling that $\psi(0)=\alpha$
we find that $\psi'(q)q_n=1$ on $\Gamma_1$.
Hence $v_n$ vanishes on $\Gamma_1$.
Since $v(x)$ satisfies the elliptic inequality (2.23), by
Hopf's second principle $v(x)$ must be a constant in $\Omega$. Because
$v(x)$ is a constant we have equality in (2.23). Then, by Remark 2.1,
$u(x)$ must be a radial function. Consequently, taking into account
conditions (3.2), $\Gamma_0$ and $\Gamma_1$ must be N-spheres.
 
Now suppose $q_0$ and $c_1$ satisfy (3.5). Let $\Gamma_1$ be the N-sphere
whose radius $r_1$ is equal to the value of $\alpha$ which solves
equation (3.20), and let $\Gamma_0$ be the N-sphere concentric with
$\Gamma_1$ whose radius $r_0$ satisfies the following equation:
$$
g(q_0^2)q_0r_0^{N-1}=r_0^N-r_1^N.                              \leqno(3.24)$$
The function $\psi(s)$ defined by
$$
g(s^2)s=(1-r_1^N\psi^{-N})\psi                                 \leqno(3.25)$$
is strictly increasing in (0,$\infty$) and satisfies the
conditions $\psi(0)=r_1$ and (by (3.24)) $\psi(q_0)=r_0$.
Its inverse function $s=\phi(r)$ satisfies equation (3.6) with $t=r$.
The theorem has been proved.
\enddemo
 
If $N=2$ and $g(q^2)=q^{p-2}$, $p>1$, then equation (3.1) is related to
the torsion problem [3,5]. In this case, conditions (3.5) read as
$$
{1\over 2}\Bigl(1-{1\over p}\Bigr)q_0^p<c_1<\Bigl(1-{1\over p}\Bigr)q_0^p.$$
Equation (3.7) becomes
$$
s^{p-1}\psi=\psi^2-\alpha^2,                                     \leqno(3.28)$$
from which we find
$$
\psi ={1\over 2}\Bigl(s^{p-1}+\bigl(s^{2p-2}+4\alpha^2\bigr)^{1\over 2}\Bigl).$$
The radius $r_1$ of the circle $\Gamma_1$ is the solution of the equation
$$
{p-1\over 2}\int_0^{q_0}\Bigl(s^{p-1}+s^{2p-2}\bigl(s^{2p-2}+
4r_1^2\bigr)^{-{1\over 2}}\Bigl)ds=c_1.$$
The radius $r_0$ of the circle $\Gamma_0$ is given by
$$
r_0={1\over 2}\Bigl(q_0^{p-1}+\bigl(q_0^{2p-2}+4r_1^2\bigr)^{1\over 2}\Bigl).$$
The solution $u(r)$ is
$$
u(r)=\int_{r_0}^r\Bigl(t-r_1^2t^{-1}\Bigr)^{1\over{p-1}}dt.$$

\enddocument 
 

\widestnumber\key{[20]}
 
\Refs
%
\ref\key 1
\by N. Garofalo, J.L. Lewis 
\paper A symmetry result related to some
overdetermined boundary value problems 
\jour Am. J. Math., 111 
\yr 1989 
\pages 9--33 
\endref
%
\ref\key 2
\by D. Gilbarg, N.S. Trudinger 
\book Elliptic partial differential
equations of second order 
\publ Springer Verlag 
\publaddr Berlin, Heidelberg, New York
\yr 1977
\endref
%
\ref\key 3
\by B. Kawohl
\paper On a family of torsional creep problems
\jour J. reine angew. Math., 410 
\yr 1990 
\pages 1--22
\endref
%
\ref\key 4
\by L.E. Payne 
\paper "Best possible" maximum principles 
\jour Math. Models and Methods in Mechanics, 
Banach Center Pubbl, 15 
\yr 1985 
\pages 609--619
\endref
%
\ref\key 5 
\by L.E. Payne, G.A. Philippin 
\paper Some applications of the maximum
principle in the problem of torsional creep 
\jour SIAM J. Appl. Math., 33 
\yr 1977
\pages 446--455
\endref
%
\ref\key 6
\by L.E. Payne, G.A. Philippin 
\paper Some maximum principles for nonlinear
elliptic equations in divergence form with applications to capillary surfaces
and to surfaces of constant mean curvature 
\jour J. Nonlinear Anal., 3 
\yr 1979 
\pages 193--211
\endref
%
\ref\key 7
\by L.E. Payne, G.A. Philippin 
\paper On maximum principles for a class
of nonlinear second order elliptic equations 
\jour J. Diff. Eq., 37 
\yr 1980 
\pages 39--48
\endref
%
\ref\key 8
\by G.A. Philippin 
\paper On a free boundary problem in electrostatics
\jour Math. Methods in Appl. Sciences, 12 
\yr 1990 
\pages 387--392
\endref
%
\ref\key 9
\by G.A. Philippin, L.E. Payne
\paper On the conformal capacity problem 
\jour Symposia Math., Vol. XXX, Academic Press 
\yr 1989 
\pages 119--136
\endref
%
\ref\key 10
\by G. Porru, F. Ragnedda 
\paper Convexity properties for solutions
of some second order elliptic semilinear equations 
\jour Applicable Analysis, 37
\yr 1990 
\pages 1--18
\endref
%
\ref\key 11 
\by M.H. Protter, H.F. Weinberger 
\book Maximum principles in differential equations 
\publ Springer-Verlag 
\publaddr Berlin, Heidelberg, New York 
\yr 1984 
\endref
%
\ref\key 12 
\by J. Serrin 
\paper A symmetry problem in potential theory  
\jour Arch. Rat. Mech. Anal., 43 
\yr 1971 
\pages 304--318 
\endref
%
\ref\key 13 
\by H.F. Weinberger 
\paper Remark on the preceding paper of Serrin
\jour Arch. Rat. Mech. Anal., 43
\yr 1971 
\pages 319--320 
\endref
\endRefs
\bye